\documentclass[reqno,b5paper]{amsart}
\usepackage{amsmath}
\usepackage{amssymb}
\usepackage{amsthm,amsxtra}
\usepackage{cite}
\usepackage{enumerate}
\usepackage[mathscr]{eucal}

\usepackage{adjustbox}
\usepackage{float}

\usepackage{tikz-cd}
\usepackage{tkz-graph}
\usepackage{tkz-berge}
\usetikzlibrary{positioning,arrows,shapes.geometric,trees}
\usepackage{graphicx}
\usepackage{mathtools}
\usetikzlibrary{positioning}
\usetikzlibrary{decorations.text}

\newtheorem{Def}{Definition}[section]
\newtheorem{Th}{Theorem}[section]
\newtheorem{Ex}{Example}[section]

\newtheorem{Lemma}{Lemma}[section]
\newtheorem{Prop}{Proposition}[section]

\newtheorem{Rem}{Remark}[section]

\newtheorem{Problem}{Problem}[section]
{}
\newcommand{\kc}{\mathcal{K}}
\newcommand{\rb}{\mathbb{R}}
\newcommand{\nb}{\mathbb{N}}

\newcommand{\hc}{\mathcal{H}}

\newcommand{\ac}{\mathcal{A}}
\newcommand{\bc}{\mathcal{B}}
\newcommand{\bk}{\mathfrak{B}}

\newcommand{\fc}{\mathcal{F}}
\newcommand{\mc}{\mathcal{M}}
\newcommand{\oc}{\mathcal{O}}

\newcommand{\tc}{\mathcal{T}}
\newcommand{\uc}{\mathcal{U}}

\newcommand{\vc}{\mathcal{V}}
\newcommand{\wc}{\mathcal{W}}
\newcommand{\yc}{\mathcal{Y}}
\newcommand{\zc}{\mathcal{Z}}

\newcommand{\obs}{\oc_{\bk^s}}
\newcommand{\bs}{{\bk^s}}
\newcommand{\gbs}{\Gamma_{\bk^s}}
\DeclareMathOperator{\Sf}{S_{fin}}
\DeclareMathOperator{\S1}{S_1}
\DeclareMathOperator{\Uf}{U_{fin}}
\DeclareMathOperator{\Gf}{G_{fin}}
\DeclareMathOperator{\G1}{G_1}
\newcommand{\bp}{\begin{proof}}
\newcommand{\ep}{\end{proof}}

\begin{document}

\title[Applications of Bornological covering properties ]{Applications of Bornological covering properties in metric spaces}

\author[D. Chandra, P. Das  and S. Das]{Debraj Chandra$^{\dag}$, Pratulananda Das$^*$  and Subhankar Das$^*$\ }

\address{\llap{$\dag$\,}Department of Mathematics, University of Gour Banga, Malda-732103, West Bengal, India}
\email{debrajchandra1986@gmail.com}

\address{\llap{*\,}Department of Mathematics, Jadavpur University, Kolkata-700032, West Bengal, India}
\email{pratulananda@yahoo.co.in, subhankarjumh70@gmail.com}

\thanks{The third author
is thankful to University Grants Commission (UGC), New Delhi-110002, India for granting UGC-NET Junior Research Fellowship (1183/(CSIR-UGC NET DEC.2017)) during the tenure of which this work was done.}


\begin{abstract}
Using the idea of strong uniform convergence\cite{bl,cmh} on bornology, Caserta, Di Maio and Ko\v{c}inac \cite{cmk} studied open covers and selection principles in the realm of metric spaces (associated with a bornology) and function spaces (w.r.t. the topology of strong uniform convergence). We primarily continue in the line initiated in \cite{cmk} and investigate the behaviour of various selection principles related to these classes of bornological covers. In the process we obtain implications among these selection principles resulting in Scheepers' like diagrams. We also introduce the notion of strong-$\bk$-Hurewicz property and investigate some of its consequences. Finally, in $C(X)$ with respect to the topology $\tau_{\bk}^s$ of strong uniform convergence, important properties like countable $T$-tightness, Reznichenko property are characterized in terms of bornological covering properties of $X$.
\end{abstract}
\maketitle
\smallskip
\noindent{\bf\keywordsname{}:} {Bornology, metric space, selection principles, open $\bs$-cover,$\gamma_\bs$-cover, $\bs$-Hurewicz property, $\bs$-Hurewicz game,$\bs$-groupability, function space, topology of strong uniform convergence.}


\section{Introduction}
We follow the notations and terminologies of \cite{arh,mcnt,engelking,hh}. We begin with the definition of the bornology on a metric space. A bornology $\bk$ on a metric space $(X,d)$ is a family of subsets of $X$ that is closed under taking finite unions, is hereditary (i.e. an ideal) and forms a cover of $X$ (see \cite{hh}).
By a base $\bk_0$ of a bornology $\bk$ we mean a subfamily of $\bk$ that is cofinal with respect to inclusion. A base is called closed (compact) if all its members are closed (compact). The smallest
bornology on $X$ is the family $\fc$ of finite subsets of $X$  while the largest bornology on $X$ is the family of all nonempty subsets of $X$.
Another important bornology is the family $\kc$ of nonempty subsets of $X$ with compact closure.\\

For two metric spaces $X$ and $Y$, $Y^X$ ($C(X, Y )$) stand for the set of all functions (continuous functions) from $X$ to $Y$. The commonly used
topologies on $C(X, Y )$ are the compact-open topology $\tau_k$, and the topology of pointwise convergence $\tau_p$. The corresponding spaces are, in general, denoted by $(C(X, Y ), \tau_k)$ (resp. $C_k(X)$ when $Y = \rb$), and $(C(X, Y ), \tau_p)$ (resp. $C_p(X)$ when $Y = \rb$).\\

In \cite{bl}, Beer and Levi had introduced the notion of strong uniform continuity on a bornology and the topology of strong uniform convergence on a bornology for function spaces. Let $(X,d)$ and $(Y ,\rho)$ be metric spaces. A mapping $f : X\rightarrow Y$ is strongly uniformly continuous on a subset $B$ of $X$ if for each $\varepsilon>0$ there is a $\delta >0$ such that $d(x_1, x_2)< \delta$ and $\{x_1, x_2\}\cap B \neq \emptyset$ imply $\rho( f(x_1),f(x_2))<\varepsilon $.
Also for a bornology $\bk$ on $X$, $f$ is called strongly uniformly continuous on $\bk$ if $f$ is strongly uniformly continuous on $B$ for each $B \in\bk$.
 Let $\bk$ be a bornology with a closed base on $X$. Then the topology of strong uniform convergence $\tau_{\bk}^s$ is determined by a uniformity on $Y^X$ with a base consisting of all sets of the form
 \[ [B,\varepsilon]^s=\{(f,g): \exists \delta> 0 \,\, \text{for every} \,\, x\in B^{\delta},(f(x),g(x))<\varepsilon \},\]
 for $B\in \bk, \varepsilon> 0$.

The topology of strong uniform convergence $\tau_{\bk}^s$ coincides with the topology of pointwise convergence $\tau_p$  if  $\bk= \fc$.\\

In another direction, \cite{cooc1} (see also \cite{cooc2}), Scheepers had began a systematic study of selection principles in topology and their relations to game theory. Study of open covers and selection principles and their inter-relationship has a long and illustrious history and readers interested in selection principles and its recent developments can consult the papers \cite{sur1,sur2,tb} where many more references can be found.\\

In \cite{cmk}, the authors, unifying the two lines of research, studied  open covers and related selection principles in the function space $C(X)(=C(X,\rb))$ with respect to the topology $\tau_{\bk}^s$ of strong uniform convergence on bornology.
In this paper we continue in that direction and use this idea of strong uniform convergence on bornology, to investigate the behaviour of various selection principles related to these classes of covers. The main objective of such investigations is to obtain ``Scheepers' like Diagrams'' which we precisely attain (Figure 1 and Figure 2) as consequences of several results where implications between these selection principles are established which is done in Section 3.
In Section 4 we introduce the notion of strong-$\bk$-Hurewicz property and strong-$\bk$-groupable cover and establish their relationships with certain selection principles and related games. In the final section, Section 5 of this paper we focus on the function space $C(X)$ endowed with the topology $\tau_{\bk}^s$. We consider several important properties like countable $T$-tightness, Reznichenko property and obtain their  characterizations in terms of bornological covering properties of $X$ among other results.

\section{Preliminaries}

Throughout the paper $(X,d)$ stands for an infinite metric space.
We first write down two classical selection principles formulated in general form in \cite{cooc1,cooc2}.
For two nonempty classes of sets $\ac$ and $\bc$ of $X$, we define\\

\noindent $\S1(\ac,\mathcal B)$: For each sequence $\{A_n:n\in \nb \}$ of elements of $\ac$, there is a sequence $\{b_n:n\in \nb\}$ such that $b_n\in A_n$ for each $n$ and $\{b_n:n\in \nb\}\in \mathcal B$.\\

\noindent $\Sf(\ac,\mathcal B)$: For each sequence $\{A_n:n\in \nb\}$ of elements of $\ac$, there is a sequence $\{B_n:n\in \nb\}$ of finite  (possibly empty) sets such that $B_n\subseteq A_n$ for each $n$ and $\bigcup_{n\in \nb}B_n\in \mathcal B$.\\

There are infinitely long games corresponding to these selection principles.\\

\noindent $\G1(\ac,\bc)$ denotes the game for two players, ONE and TWO, who play a round for each positive integer $n$. In the $n$-th round ONE chooses a set $A_n$ from $\ac$ and TWO responds by choosing an element $b_n\in  A_n$. TWO wins the play $\{A_1,b_1, \dotsc, A_n, b_n, \dotsc \}$ if $\{b_n :n\in \nb\}\in \bc$. Otherwise ONE wins.\\

\noindent $\Gf(\ac,\bc)$ denotes the game where in the $n$-th round ONE chooses a set $A_n$ from $\ac$ and TWO responds by choosing a finite (possibly empty) set $B_n\subseteq A_n$. TWO wins the play $\{A_1,B_1, \dotsc,A_n,B_n, \dotsc  \}$ if $\bigcup_{n\in \nb}B_n\in \bc$. Otherwise ONE wins.\\

We also consider the following selection principles.\\

\noindent $\Uf(\ac,\mathcal B)$: For each sequence $\{A_n:n\in \nb\}$ of elements of $\ac$, there is a sequence $\{B_n:n\in \nb\}$ of finite  (possibly empty) sets such that $B_n\subseteq A_n$ for each $n$ and either $\{\cup B_n:n\in \nb\}\in \mathcal B$ or for some $n$, $\cup B_n=X$ (from \cite{cooc1,cooc2}).\\

\noindent $\binom{\ac}{\bc}$: For each element $A$ of $\ac$, there is a set $B$ such that $B\subseteq A$ and $B\in \bc$ (from \cite{tb04}).\\

We now recall some classes of open covers.
Let $\oc$ denote the collection of all open covers of $X$. An open cover $\uc$ of $X$ is called a $\gamma$-cover \cite{cooc1} (see also \cite{jgzsn}) if $\uc$ is infinite and each $x\in X$ belongs to all but finitely many members of $\uc$. The collection of all $\gamma$-covers of $X$ is denoted by $\Gamma$.
An open cover $\uc$ of $X$ is said to be groupable \cite{coocvii} if it can be represented as a union of countably many finite pairwise disjoint sets $\uc_n$ such that each $x$ in $X$ belongs to $\cup \uc_n$ for all but finitely many $n\in \nb$. The collection of all groupable open covers is denoted by $\oc^{gp}$.
For $x$ in $X$, we denote $\Omega_x=\{A\subseteq X:x\in \overline{A}\setminus A\}$ \cite{coocvii}.
A countable set $A\in \Omega_x$ is groupable \cite{coocvii} if there is a partition $\{B_n:n\in \nb\}$ of $A$ into pairwise disjoint finite sets such that each neighbourhood of $x$ has non-empty intersection with all but finitely many $B_n$. The family of all groupable elements of $\Omega_x$ is denoted by $\Omega_x^{gp}$. For each $x$ in $X$, $\Sigma_x$ denotes the collection of all sequences that converges to $x$ \cite{kock2}.
Throughout the paper we assume that all the open covers of $X$ are countably infinite.\\

Next we recall certain notions which we would deal with in the final section.
 $X$ is said to have countable fan tightness at $x\in X$ \cite{arh86} if $X$ satisfies $\Sf(\Omega_x,\Omega_x)$.
 $X$ is said to have countable $T$-tightness \cite{j}  if for each uncountable regular cardinal $\rho$ and each increasing sequence $\{A_\alpha:\alpha<\rho\}$ of closed subsets of $X$ the set $\cup\{A_\alpha:\alpha<\rho\}$ is closed.
For $x\in X$, a family $\ac$ of subsets of $X$ is called a $\pi$-network at $x$ if every neighbourhood of $x$ contains some element of $\ac$. $X$ is called a Pytkeev space (\!\cite{pyt}, see \cite{mt,s} for more details) if $x\in \overline{A}\setminus A$ and $A\subset X$ imply the existence of a countable $\pi$-network at $x$ of infinite subsets of $A$. $X$ is called  weakly Fr\'{e}chet-Urysohn \cite{fed} if for $x\in \overline{A}\setminus A$ and $A\subset X$ there exists a countable infinite disjoint family of finite subsets $\{A_n:n\in \nb\}$ of $A$ such that every neighbourhood of $x$ intersects $A_n$ for all but finitely many $n\in \nb$. This property is also known as Reznichenko property of $X$.\\

Let $\bk$ be a bornology for a metric space $(X,d)$ with closed base. For $B\in\bk$ and $\delta>0$, let $B^{\delta}=\bigcup_{x\in B} S(x,\delta )$, where $S(x,\delta)=\{y\in X:d(x,y)<\delta\}$. It is easy to verify that $\overline{B}\in \bk$ for every $B\in \bk$ and $\overline{B^{\delta}}\subseteq B^{2\delta}$ for every $B\in\bk$ and $\delta>0$.
A cover $\uc$ of $X$ is said to be a strong $\bk$-cover (in short, $\bs$-cover) of $X$ \cite{cmh} if $X$ is not in $\uc$ and for each $B\in \bk$ there exist $U\in \uc$ and $\delta>0$ such that $B^\delta\subseteq U$. A $\bs$-cover $\uc$ is said to be an open $\bs$-cover if the members of $\uc$ are open sets.
The collection of all open $\bs$-covers of $X$ is denoted by $\oc_{\bk^s}$. Clearly $\oc_{\bk^s}\subseteq \oc$.
An open cover $\uc=\{U_n:n\in \nb\}$ of $X$ is said to be a $\gamma_{\bk^s}$-cover \cite{cmk}(see also \cite{cmh}) if it is infinite and for each $B\in \bk$ there is a $n_0\in  \nb$ and a sequence $\{\delta_n:n\ge  n_0\}$ of positive real numbers
such that $B^{\delta_n}\subseteq U_n$ for all $n\ge  n_0$.
The collection of all $\gamma_{\bk^s}$-covers of $X$ is denoted by $\Gamma_{\bk^s}$.
Throughout the paper we use the convention that if $\bk$ is a bornology on a metric space $X$, then $X\notin \bk$.
\section{Selection Principles and Bornological covering properties}

\subsection{Certain observations on covers}

We start this section with a basic characterization of $\bk^s$-covers of $X$.
\begin{Prop}
\label{P1}
Let $\bk$ be a bornology on $X$ with closed base. For an open cover $\uc=\{U_n:n\in \nb\}$ of $X$, the following conditions are equivalent.\\
\noindent $(1)$ $\uc$ is an open $\bk^s$-cover of $X$.\\
\noindent $(2)$ For every $B\in \bk$ there are $\delta_n>0$ and $U_n\in \uc$ such that $B^{\delta_n}\subseteq U_n$ for infinitely many $n\in \nb$.
\end{Prop}

\begin{proof} We only give the proof of $(1)\Rightarrow (2)$.
For $B\in \bk$, choose $x_1\in X\setminus B$ (as $X\notin \bk$). Then there is a $B_{x_1}\in\bk$ such that $x_1\in B_{x_1}$ and $B\subsetneqq B\cup B_{x_1}\in \bk$. Consequently one can find  $\delta_{n_1}>0$ and $U_{n_1}\in \uc$ such that $(B\cup B_{x_1})^{\delta_{n_1}}\subseteq U_{n_1}$, i.e. $B^{\delta_{n_1}}\subseteq U_{n_1}$. Let $x_2\in X\setminus(B\cup B_{x_1})$. Then there is a $B_{x_2}\in \bk$ such that $x_2\in B_{x_2}$ and $B\cup B_{x_1}\subsetneqq (B\cup B_{x_1}\cup B_{x_2})\in \bk$. Which consequently assures the existence of $\delta_{n_2}>0$ and $U_{n_2}\in \uc$ such that $(B\cup B_{x_1}\cup B_{x_2})^{\delta_{n_2}}\subseteq U_{n_2}$, i.e. $B^{\delta_{n_2}}\subseteq U_{n_2}$. Continuing in this process, we obtain an increasing sequence $n_1<n_2< \cdots$ with positive real numbers $\delta_{n_k}$ and $U_{n_k}\in \uc$ such that $B^{\delta_{n_k}}\subseteq U_{n_k}$ for all $k$. Thus $(2)$ holds.
\end{proof}

Another observation about the $\bs$-cover is the following.
\begin{Lemma}
\label{L2}
Let $\bk$ be a bornology on $X$ with closed base. Let $\{\uc_n:n\in \nb\}$ be a sequence of open ${\bk^s}$-covers of $X$ and  $\uc_n=\{U^n_j:j\in\nb\}$ for each $n$.
If $\vc_n=\{ U_{k_1}^1\cap U_{k_2}^2\cap \dotsc \cap U_{k_n}^n:U_{k_i}^i\in \uc_i, k_i\in \nb,i=1,2, \dotsc, n\}$, then for each $n\in \nb$, $\vc_n$ is an open ${\bk^s}$-cover of $X$.
\end{Lemma}

\begin{proof}
Let $\uc_n=\{U_k^n:k\in \nb\}$. For $B\in \bk$ there exist  $\delta_{k_i}^i>0$ and $U_{k_i}^i\in \uc_i$ such that $B^{\delta_{k_i}^i}\subseteq U_{k_i}^i$ for each $i=1,2, \dotsc, n$.
Choose $\delta_k=\min\{\delta_{k_i}^i:i=1,2, \dotsc, n\}$. Then we have $B^{\delta_k}\subseteq U_{k_1}^1\cap U_{k_2}^2\cap \dotsc \cap U_{k_n}^n\in \vc_n$. Hence $\vc_n$ is an open $\bk^s$-cover of $X$ for each $n\in \nb$.
\end{proof}

Coming to $\gamma_{\bk^s}$-covers, it is easy to observe that every $\gamma_{\bk^s}$-cover is a $\gamma$-cover. The following example shows that the class of $\gamma_\bs$-covers is properly contained in the class of $\gamma$-covers.

\begin{Ex}
\label{E1}
Consider $X=\mathbb{R}$ with the Euclidean metric $d$ and the bornology $\bk=\{(-\infty,x):x\in \mathbb{R}\}$. Let $\uc=\{(-m,m):m\in \mathbb{N}\}$. Then $\uc$ is a $\gamma$-cover of $X$.
If $B=(-\infty,x_0)\in \bk$ and $\delta >0$, then $B^\delta=(-\infty, x_0+\delta)$. Clearly $B^\delta\nsubseteq (-m,m)$ for any $m\in \mathbb{N}$ and any $\delta >0$. Hence $\uc$ can not be a $\gamma_{\bk^s}$-cover of $X$.
\end{Ex}

Thus the following  inclusion relations among covers can be observed.
$\Gamma_{\bk^s}\subsetneqq \oc_{\bk^s}\subsetneqq \oc, \,\Gamma_{\bk^s}\subsetneqq \Gamma$.

A key observation about $\gamma_\bs$-covers is the following in the line of \cite{cmk} which will come to our use.
\begin{Lemma}
Let $\bk$ be a bornology on $X$ with closed base. If a $\gamma_\bs$-cover of $X$ is partitioned into finitely many pieces, then each of them is again a $\gamma_\bs$-cover.
\end{Lemma}


\begin{Lemma}[cf. \cite{cmk}]
\label{L1}
 Let $\bk$ be a bornology on $X$ with closed base.\\
\noindent $(1)$ Any infinite subset of a $\gamma_{\bk^s}$-cover is a $\gamma_{\bk^s}$-cover.\\
\noindent $(2)$  Let $\{\uc_n:n\in \nb\}$ be a sequence of $\gamma_{\bk^s}$-covers of $X$ and  $\uc_n=\{U^n_j:j\in\nb\}$ for each $n$. Then $\vc_n=\{U_k^1\cap U_k^2\cap \dotsc \cap U_k^n:k\in \nb, U_k^i\in \uc_i,  i=1,2, \dotsc, n\}$ is a $\gamma_{\bk^s}$-cover of $X$.
\end{Lemma}

Finally we observe that one can construct a $\gamma_\bs$-cover from an open $\bs$-cover in the following way.
\begin{Lemma}
\label{L3}
Let $\bk$ be a bornology on $X$ with closed base. Let $\uc=\{U_n:n\in \nb\}$ be an open ${\bk^s}$-cover of $X$. If $\vc=\{V_n:n\in \nb\}$ where $V_n=\cup_{i=1}^{n}U_i$, then $\vc$ is a $\gamma_{\bk^s}$-cover of $X$.
\end{Lemma}

\begin{proof}
 Let $B\in \bk$. Since $\uc$ is an open $\bk^s$-cover of $X$, there exist $U_{k_0}\in \uc$ and $\delta>0$ such that $B^\delta\subseteq U_{k_0}$. Then $B^\delta\subseteq \cup_{i=1}^{k}U_i=V_k$ for all $k\ge k_0$. If we take  $\delta_n=\delta$ for $n\ge k_0$, then $B^{\delta_n}\subseteq V_n$ for all $n\ge k_0$. Hence $\vc$ is a $\gamma_{\bk^s}$-cover of $X$.
\end{proof}

\subsection{Implications among selection principles and Scheepers' like diagrams}

Taking $\ac,\bc\in\{\oc,\Gamma,\obs,\gbs\}$ and $\prod\in\{\S1,\Sf,\Uf\}$, we establish the equivalences among the selection principles $\prod(\ac,\bc)$ in the next few results.
\begin{Th}
\label{T1}
Let $\bk$ be a bornology on $X$ with closed base. The following statements hold.\\
\noindent$(1)$ $\S1(\Gamma_{\bk^s},\Gamma_{\bk^s})=\Sf(\Gamma_{\bk^s},\Gamma_{\bk^s})$\\
\noindent$(2)$ $\S1(\Gamma_{\bk^s},\Gamma)=\Sf(\Gamma_{\bk^s},\Gamma)$\\
\noindent$(3)$ $\S1(\Gamma,\Gamma_{\bk^s})=\Sf(\Gamma,\Gamma_{\bk^s})$.
\end{Th}

\begin{proof}
 We  only prove $(1)$  and omit the remaining proofs as they are analogous.\\

Clearly $\S1(\Gamma_{\bk^s},\Gamma_{\bk^s})$ implies $\Sf(\Gamma_{\bk^s},\Gamma_{\bk^s})$.

Let $X$ satisfy $\Sf(\Gamma_{\bk^s},\Gamma_{\bk^s})$ and let $\{\uc_n:n\in \nb\}$ be a sequence of $\gamma_{\bk^s}$-covers of $X$. For each $n\in \nb$, we enumerate $\uc_n$ bijectively as $\uc_n=\{U_k^n:n\in \nb\}$. Define $\vc_n=\{U_k^1\cap U_k^2\cap \dotsc \cap U_k^n:k\in \nb\}$ for each $n$. By Lemma \ref{L1}(2), $\vc_n$ is a $\gamma_{\bk^s}$-cover of $X$. Now applying $\Sf(\Gamma_{\bk^s},\Gamma_{\bk^s})$ to this sequence $\{\vc_n:n\in \nb\}$ of $\gamma_{\bk^s}$-covers of $X$, we can find a sequence $\{\wc_n:n\in \nb\}$ of finite subsets with $\wc_n\subseteq \vc_n$ such that $\cup_{n=1}^\infty\wc_n$  is a $\gamma_{\bk^s}$-cover of $X$. Thus it is possible to find a sequence $n_1<n_2< \cdots $ of positive integers such that for each $j$, $\wc_{n_j}\setminus \cup_{i<j}\wc_{n_i}$ is non-empty and their union $\cup_{k\in \nb}\wc_{n_k}$ is a $\gamma_{\bk^s}$-cover of $X$ by Lemma \ref{L1}(1).
Now for each $j$, choose $m_j$ such that $V_{m_j}^{n_j}$ is an element of $\wc_{n_j}\setminus \cup_{i<j}\wc_{n_i}$. Then $\{V_{m_j}^{n_j}:j=1,2, \dotsc \}$ is a $\gamma_{\bk^s}$-cover of $X$  by Lemma \ref{L1}(1). Define $U_n=U_{m_{k+1}}^n$ for each $n\in (n_k,n_{k+1}]$. Consider the set $\{U_n:n=1,2, \dotsc \}$.
For $B\in \bk$ there exist $k_0\in \nb$ and a sequence $\{\delta_{n_k}:k\ge k_0\}$ of positive reals such that $B^{\delta_{n_{k+1}}}\subseteq V_{m_{k+1}}^{n_{k+1}}$ for all  $k\ge  k_0$, i.e.
$B^{\delta_{n_{k+1}}}\subseteq U_{m_{k+1}}^1\cup U_{m_{k+1}}^2\cup \cdots \cup U_{m_{k+1}}^{n_{k+1}}$. For $n\in (n_k,n_{k+1}]$, define $\delta_n= \delta_{n_{k+1}}$. Then $B^{\delta_n}\subseteq U_n$ for all $n\ge  n_{k_0}$ and $k\ge k_0$. So $\{U_n:n\in \nb\}$ is a $\gamma_{\bk^s}$-cover of $X$. Hence $X$ satisfies $\S1(\Gamma_{\bk^s},\Gamma_{\bk^s})$.
\end{proof}

\begin{Th}
\label{T2}
Let $\bk$ be a bornology on $X$ with closed base. The following statements hold.\\
\noindent$(1)$ $\S1(\oc_{\bk^s},\Gamma_{\bk^s})=\Sf(\oc_{\bk^s},\Gamma_{\bk^s})$\\
\noindent$(2)$ $\S1(\oc,\Gamma_{\bk^s})=\Sf(\oc,\Gamma_{\bk^s})$\\
\noindent$(3)$ $\S1(\oc_{\bk^s},\Gamma)=\Sf(\oc_{\bk^s},\Gamma)$.
\end{Th}

\begin{proof}
We only present an outline of the proof of $(1)$.
Let $X$ satisfy  $\Sf(\oc_{\bk^s},\Gamma_{\bk^s})$ and let $\{\uc_n:n\in \nb\}$ be a sequence of open ${\bk^s}$-covers of $X$, where $\uc_n=\{U_k^n:k\in \nb\}$.
Now define $\vc_n=\{U_{k_1}^1\cap U_{k_2}^2\cap \dotsc \cap U_{k_n}^n:U_{k_i}^i\in \uc_i, k_i\in \nb,i=1,2, \dotsc, n\}$. By Lemma \ref{L2}, for each $n\in \nb$, $\vc_n$ is an open ${\bk^s}$-cover of $X$. Now applying $\Sf(\oc_{\bk^s},\Gamma_{\bk^s})$ to the sequence $\{\vc_n:n\in \nb\}$ and proceeding as in the proof of  Theorem \ref{T1}$(1)$, we can conclude that $X$ satisfies $\S1(\oc_{\bk^s},\Gamma_{\bk^s})$.
\end{proof}

\begin{Th}
\label{T3}
Let $\bk$ a bornology on $X$ with closed base. The following statements hold.\\
\noindent$(1)$ $\Sf (\Gamma_{\bk^s},\oc)=\Uf (\oc_{\bk^s},\oc)$\\
\noindent$(2)$ $\Sf(\oc_{\bk^s},\oc)=\Sf(\Gamma_{\bk^s},\oc)$\\
\noindent$(3)$ $\Sf(\oc_{\bk^s},\oc)=\Uf(\Gamma_{\bk^s},\oc)$.
\end{Th}

\bp
 We will again give proof of only $(1)$.
Let $X$ satisfy $\Sf(\Gamma_{\bk^s},\oc)$ and let $\{\uc_n:n\in \nb\}$ be a sequence of open $\bk^s$-covers of $X$. Let $\uc_n=\{U_k^n:k\in \nb\}$. Now consider the collection $\vc_n=\{V_k^n:k\in \nb\}$ where $V_k^n=U_1^n\cup U_2^n\cup \cdots \cup U_k^n$. Then by Lemma \ref{L3}, $\{\vc_n:n\in \nb\}$ is a sequence of $\gamma_{\bk^s}$-covers of $X$. Applying $\Sf(\Gamma_{\bk^s},\oc)$ to this sequence, we get a finite subset $\wc_n\subseteq \vc_n$ such that $\cup_{n=1}^{\infty}\wc_n$ is an open cover of $X$. Now deconstructing the members of each $\wc_n$, we find for each $n$, a finite subset $\zc_n\subseteq \uc_n$ such that $\{\cup \zc_n:n\in \nb\}$ is an open cover of $X$. For $x\in X$ there exists $V_{k_0}^{n_0}\in \cup_{n=1}^{\infty}\wc_n$ such that $x\in V_{k_0}^{n_0}=U_1^{n_0}\cup \cdots \cup U_{k_0}^{n_0}$, i.e. $x\in U_i^{n_0}$ for some $i\in \{1,2, \dotsc, k_0\}$, i.e. $x\in U_i^{n_0}\subseteq \cup\zc_{n_0}$ (as $U_i^{n_0}\in \zc_{n_0}$). Hence $\{\cup\zc_{n_0}:n\in \nb\}$ is an open cover of $X$. So $X$ satisfies $\Uf(\oc_{\bk^s},\oc)$.

Conversely let $X$ satisfy $\Uf(\oc_{\bk^s},\oc)$ and $\{\uc_n:n\in \nb\}$ be a sequence of $\gamma_{\bk^s}$-covers of $X$. Applying $\Uf(\oc_{\bk^s},\oc)$ to $\{\uc_n:n\in \nb\}$, we obtain for each $n\in \nb$, a finite set $\vc_n\subseteq \uc_n$  such that $\{\cup \vc_n:n\in \nb\}$ is an open cover of $X$. For $x\in X$, $x\in \cup\vc_{k_0}$ for some $k_0\in \nb$, i.e. $x\in U\subseteq \cup \vc_{k_0}$ for some $U\in \vc_{k_0}\subseteq \uc_{k_0}$. So $\cup_{n=1}^{\infty}\vc_k$ is an open cover of $X$. Hence $X$ satisfies $\Sf(\Gamma_{\bk^s},\oc)$.
\ep

\begin{Th}
\label{T4}
Let $\bk$ be a bornology on $X$ with closed base. The following statements hold.\\
\noindent$(1)$ $\Uf(\oc_{\bk^s},\oc_{\bk^s})=\Uf(\Gamma_{\bk^s},\oc_{\bk^s})$.\\
\noindent$(2)$ $\Uf(\oc_{\bk^s},\Gamma_{\bk^s})=\Uf(\Gamma_{\bk^s},\Gamma_{\bk^s})$.\\
\noindent$(3)$ $\Uf(\oc_{\bk^s},\oc)=\Uf(\Gamma_{\bk^s},\oc)$ \\
\noindent$(4)$ $\Uf(\oc,\oc_{\bk^s})=\Uf(\Gamma,\oc_{\bk^s})$\\
\noindent$(5)$ $\Uf(\oc_{\bk^s},\Gamma)=\Uf(\Gamma_{\bk^s},\Gamma)$\\
\noindent$(6)$ $\Uf(\oc,\Gamma_{\bk^s})=\Uf(\Gamma,\Gamma_{\bk^s})$.
\end{Th}

\bp
$(1)$.
If $X$ satisfies $\Uf(\oc_{\bk^s},\oc_{\bk^s})$, then clearly $X$ also satisfies  $\Uf(\Gamma_{\bk^s},\oc_{\bk^s})$.

Now let $X$ satisfy $\Uf(\Gamma_{\bk^s},\oc_{\bk^s})$ and $\{{\uc_n}:n\in \nb\}$ be a sequence of open ${\bk^s}$-covers of $X$. For each $n$, we enumerate $\uc_n$ as  $\uc_n=\{U_k^n:k=1,2, \dotsc \}$. Consider the collection $\vc_n$ whose $m$-th member is $V_m=\cup_{i=1}^{m}U_i^n$. By Lemma \ref{L3}, $\vc_n$ is a $\gamma_{\bk^s}$-cover of $X$. So $\{\vc_n:n\in \nb\}$ is a sequence of $\gamma_{\bk^s}$-covers of $X$. Applying $\Uf(\Gamma_{\bk^s},\oc_{\bk^s})$ to this sequence, we get a finite set $\wc_n\subseteq \vc_n$  such that $\{\cup \wc_n:n\in \nb\}$ is an open ${\bk^s}$-cover of $X$. Now deconstructing the members of each $\wc_n$, we find, for each $n$, a finite subset $\tc_n$ of $\uc_n$ with $\cup \wc_n=\cup \tc_n$.
Since $\{\cup\wc_n:n\in \nb\}$ is an open $\bk^s$-cover of $X$ so $\{\cup \tc_n:n\in \nb\}$ is also an open $\bk^s$-cover of $X$. Hence $X$ satisfies $\Uf(\oc_{\bk^s},\oc_{\bk^s})$.\\
The remaining proofs are analogous and so are omitted.
\ep

Combining Theorem \ref{T1}, Theorem \ref{T2}, Theorem \ref{T3} and Theorem \ref{T4}, the following implication diagram (Figure 1) can be easily verified.

\begin{figure}[h]
\begin{adjustbox}{max width=\textwidth,max height=\textheight,keepaspectratio,center}
\begin{tikzcd}[column sep=1.5ex,row sep=10ex,arrows={crossing over}]
%
&&&
\Uf(\Gamma,\obs)&&&&&
\\
&&&
\Uf(\oc,\obs) \arrow[u,thick,equals]\arrow[ddl,leftarrow] &&&&&
\\
&&\Uf(\Gamma,\gbs) \arrow[uur] \arrow[rr] &&
\Uf(\gbs,\gbs) \arrow[rr] \arrow[ddl,leftarrow] &&
\Uf(\gbs,\Gamma) \arrow[ddl,leftarrow]\arrow[rr]&&
\Uf(\gbs,\oc)\arrow[ddl,thick,equals]
\\
&&\Uf(\oc,\gbs) \arrow[ddl,leftarrow] \arrow[u,thick,equals] \arrow[rr]  &&
\Uf(\obs,\gbs) \arrow[u,thick,equals] \arrow[rr] \arrow[ddl,leftarrow] &&
\Uf(\obs,\Gamma) \arrow[u,thick,equals]\arrow[rr]\arrow[ddl,leftarrow] &&
\Uf(\obs,\oc)\arrow[u,thick,equals]
\\
&\Sf(\Gamma,\gbs) \arrow[uur]\arrow[rr] &&
\Sf(\gbs,\gbs) \arrow[ddl,thick,equals] \arrow[rr] &&
\Sf(\gbs,\Gamma) \arrow[ddl,thick,equals] \arrow[rr] &&
\Sf(\gbs,\oc) \arrow[ddl,leftarrow]&
\\
&
\Sf(\oc,\gbs) \arrow[u] \arrow[rr]  \arrow[ddl,thick,equals]&&
\Sf(\obs,\gbs) \arrow[u] \arrow[rr] \arrow[ddl,thick,equals]&&
\Sf(\obs,\Gamma) \arrow[u] \arrow[rr]\arrow[ddl,thick,equals]&&
\Sf(\obs,\oc) \arrow[u,thick,equals]\arrow[uur]\arrow[ddl,leftarrow]&
\\
\S1(\Gamma,\gbs)\arrow[uur,thick,equals] \arrow[rr] &&
\S1(\gbs,\gbs) \arrow[rr]&&
\S1(\gbs,\Gamma)\arrow[rr]&&
\S1(\gbs,\oc)&&
\\
 \S1(\oc,\gbs)\arrow[u] \arrow[rr] &&
\S1(\obs,\gbs) \arrow[u]\arrow[rr]&&
\S1(\obs,\Gamma)\arrow[u]\arrow[rr]&&
\S1(\obs,\oc)\arrow[u]\arrow[uur]&&
\end{tikzcd}
\end{adjustbox}
 \caption{}
\end{figure}

Furthermore if one considers only the classes $\obs$ and $\gbs$, then one can obtain the following diagrammatic representation (Figure 2).
\newpage

\begin{figure}[h]
\begin{adjustbox}{width=\textwidth,height=\textheight,keepaspectratio,center}
\begin{tikzcd}[column sep=8ex,row sep=8ex,arrows={crossing over}]
%
&&
\Uf(\gbs,\gbs)\arrow[rr]&&
\Uf(\gbs,\obs)\arrow[ddl,leftarrow]
\\
&&
\Uf(\obs,\gbs) \arrow[u,thick,equals]\arrow[ddl,leftarrow]\arrow[rr] &&
\Uf(\obs,\obs)\arrow[u,thick,equals]
\\
&\Sf(\gbs,\gbs) \arrow[uur] \arrow[rr] &&
\Sf(\gbs,\obs)\arrow[ddl,leftarrow] &&
\\
&\Sf(\obs,\gbs)\arrow[ddl,thick,equals] \arrow[u] \arrow[rr] &&
\Sf(\obs,\obs)\arrow[u]\arrow[uur] &&
\\
\S1(\gbs,\gbs)\arrow[uur,thick,equals]\arrow[rr]&&
\S1(\gbs,\obs)&&
\\
\S1(\obs,\gbs) \arrow[u]\arrow[rr] &&
\S1(\obs,\obs)\arrow[u]\arrow[uur]&&
\end{tikzcd}
\end{adjustbox}
 \caption{}
\end{figure}

\subsection{Game theoretic characterizations of some selection principles}

 The study of topological properties and their relations to Ramsey theory, function spaces, and other related topics can be described in terms of  topological games. Rich surveys are available in \cite{t74,t75,t84,paw,cooc1,cooc2} and reference therein. Numerous properties in selection principles are interconnected with two player game and these properties are equivalent to the fact that the first player does not have a wining strategy in that game. From time to time these game theoretic observations become useful tools to derive results related to those properties in selection principles.

We first present the following game theoretic characterization of $\S1(\gbs,\gbs)$ which is motivated from the similar characterization of $\S1(\Gamma,\Gamma)$ (\!\cite[Theorem 26]{cooc1}). For the convenience of readers, we completely outline the proof here.
\begin{Th}
\label{T5}
Let $\bk$ a bornology on $X$ with closed base. The following conditions are equivalent.\\
\noindent$(1)$ $X$ satisfies $\S1(\Gamma_{\bk^s},\Gamma_{\bk^s})$.\\
\noindent$(2)$ ONE has no wining strategy in the game $\G1 (\Gamma_{\bk^s},\Gamma_{\bk^s})$.
\end{Th}

\bp

\noindent$(1)\Rightarrow (2)$.
Suppose that $F$ is a strategy for ONE. Let the first move of ONE be  $F(X)=\uc_1=\{U_{(n)}:n\in \nb\}$, a $\gamma_{\bk^s}$-cover of $X$. TWO responds by selecting $U_{(n_1)}\in \uc_1$. Then the 2nd move of  ONE is $F(U_{(n_1)})\setminus \{U_{(n_1)}\}=\{U_{(n_1,m)}:m\in\nb\}\in  \Gamma_{\bk^s}$.
TWO responds by selecting $U_{(n_1,n_2)}\in F(U_{(n_1)})\setminus \{U_{(n_1)}\}$. Let the $(k+1)$-th move of ONE be $F(U_{(n_1)}, \dotsc, U_{(n_1, \dotsc ,n_k)})\setminus \{U_{(n_1)}, \dotsc, U_{(n_1, \dotsc ,n_k)}\}=\{U_{(n_1, \dotsc ,n_k,m)}:m\in \nb\}\in \Gamma_{\bk^s}$.
Hence for each finite sequence $\sigma$, the collection $\{U_{\sigma\frown (m)}:m\in \nb\}$ is a $\gamma_{\bk^s}$-cover of $X$. Since $X$ satisfies $S_1(\Gamma_{\bk^s},\Gamma_{\bk^s})$, for each $\sigma$ there exist a $n_{\sigma}\in \nb$ such that $\{U_{\sigma\frown (n_\sigma)}:\sigma \quad \text{is a finite sequence}\}$ is a $\gamma_{\bk^s}$-cover.
Define inductively $n_1,n_2,\dotsc$ so that $n_1=n_\emptyset$ and $n_{k+1}=n_{(n_1, \dotsc ,n_k)}$.
Then the sequence $U_{(n_1)},U_{(n_1,n_2)}, \dotsc ,U_{(n_1, \dotsc ,n_k)}, \dotsc$ of moves of TWO is a $\gamma_{\bk^s}$-cover by Lemma \ref{L1}. Thus $F$ is not a winning strategy in the game $\G1 (\Gamma_{\bk^s},\Gamma_{\bk^s})$.\\

\noindent$(2)\Rightarrow (1)$.
Suppose that $X$ does not satisfy $S_1(\Gamma_{\bk^s},\Gamma_{\bk^s})$. Then there exists a sequence  $\{\uc_n :n\in \nb\}$ of $\gamma_{\bk^s}$-covers of $X$ such that for any choice of $U_n\in \uc_n$, $\{U_n:n\in \nb\}$ is not a $\gamma_{\bk^s}$-cover of $X$.
Now we define a strategy $F$ for ONE as follows. Let the first move of ONE be  $\uc_1$. TWO responds by selecting $U_1\in \uc_1$. Let the second move of ONE be  $\uc_2$. TWO responds by selecting $U_2\in \uc_2$. At the $n$-th stage ONE's move is $\uc_n$ and TWO selects $U_n\in \uc_n$. Since $F$ is not a winning strategy for ONE, the collection $\{U_n:n\in \nb\}$ becomes a $\gamma_{\bk^s}$-cover, which is a contradiction. Hence $(1)$ is satisfied.
\ep

In the following a similar characterization can be observed.
\begin{Th}
\label{T6}
Let $\bk$ a bornology on $X$ with closed base. The following conditions are equivalent.\\
\noindent$(1)$  $X$ satisfies $\Sf(\Gamma_{\bk^s},\Gamma_{\bk^s})$.\\
\noindent$(2)$  ONE has no wining strategy in the game $\Gf(\Gamma_{\bk^s},\Gamma_{\bk^s})$.
\end{Th}

Using Proposition \ref{P1} and replacing $\gamma_\bs$-covers by open $\bs$-covers in the first component of the selection principle in Theorem \ref{T5}, a similar characterization can be obtained.
\begin{Th}
\label{T7}
Let $\bk$ be a bornology on $X$ with closed base. The following conditions are equivalent.\\
\noindent$(1)$ $X$ satisfies $\S1(\oc_{\bk^s},\Gamma_{\bk^s})$.\\
\noindent$(2)$  ONE has no wining strategy in the game $\G1 (\oc_{\bk^s},\Gamma_{\bk^s})$.
\end{Th}

\section{The strong $\bk$-Hurewicz property}

We first introduce the definition of strong $\bk$-Hurewicz property.
\begin{Def}
\label{D1}
Let $\bk$ be a bornology on $X$ with closed base. $X$ is said to have strong $\bk$-Hurewicz property (in short, $\bk^s$-Hurewicz property) if for each sequence $\{\uc_n:n\in \nb\}$ of open ${\bk^s}$-covers of $X$, there is a sequence  $\{\vc_n:n\in \nb\}$ where $\vc_n$ is a finite subset of $\uc_n$ for each $n\in \nb$, such that for every $B\in \bk$ there exist a $n_0\in \nb$ and a sequence $\{\delta_n:n\ge  n_0\}$ of positive real numbers satisfying $B^{\delta_n}\subseteq U$ for some $U\in \vc_n$ for all $n\ge  n_0$.
\end{Def}

\begin{Def}
The strong $\bk$-Hurewicz game (in short, $\bk^s$-Hurewicz game) on $X$ is defined as follows.
Two players named ONE and TWO play an infinite long game. In the $n$-th inning ONE selects an open $\bk^s$-cover $\uc_n$ of $X$, TWO responds by choosing a finite set $\vc_n\subseteq \uc_n$. TWO wins the play: $\uc_1, \vc_1, \uc_2, \vc_2, \dotsc ,\uc_n, \vc_n, \dotsc$ if for each $B\in \bk$ there exist a $n_0\in \nb$ and a sequence $\{\delta_n:n\ge  n_0\}$ of positive real numbers satisfying $B^{\delta_n}\subseteq U$ for some $U\in \vc_n$ for all $n\ge  n_0$. Otherwise ONE wins.
\end{Def}

The $\bs$-Hurewicz property can be completely characterized by the $\bs$-Hurewicz game.
\begin{Th}
\label{T8}
Let $\bk$ be a bornology on $X$ with closed base. The following conditions are equivalent.\\
\noindent$(1)$ $X$ has the $\bk^s$-Hurewicz property.\\
\noindent$(2)$ ONE has no wining strategy in the $\bk^s$-Hurewicz game on $X$.
\end{Th}

\bp \noindent$(1)\Rightarrow (2)$.
Let $X$ have the $\bk^s$-Hurewicz property. Let $F$ be the strategy for ONE in the $\bs$-Hurewicz game on $X$. The first move of ONE is $F(X)=\{U_{(n)}:n\in \nb\}$, which is an open-$\bk^s$-cover of $X$. TWO responds by selecting a finite set $\vc_{(n_1)}=\{U_{(n)}:n\le n_1\}$. Note that if we discard finitely many elements from an open $\bs$-cover, the remaining members still form an open $\bs$-cover by Proposition \ref{P1}. Suppose that for each finite sequence $\sigma$ of positive integers of length at most $m$, $U_\sigma$ has been defined. Define $F\left(\vc_{(n_1)}, \dotsc ,\vc_{(n_1, \dotsc ,n_{m-1})}\right)\setminus \{\vc_{(n_1)}\cup \cdots \cup\vc_{(n_1, \dotsc ,n_{m-1})}\}=\{U_{(n_1, \dotsc ,n_{m-1},k)}:k\in \nb\}$. Clearly for each finite sequence $\sigma$ of positive integers, $\uc_\sigma=\{U_{\sigma\frown (n)}:n\in \nb\}$ is an open $\bs$-cover of $X$. By $(1)$, we can find for each $\sigma$, a $n_\sigma\in \nb$ and a finite set $\vc_\sigma=\{U_{\sigma\frown (n_{\sigma})}:n\le n_\sigma\}$ such that $\{\vc_\sigma:\sigma \text{ is a finite sequence}\}$  witnesses the $\bs$-Hurewicz property of $X$.

Define inductively a sequence $n_1=n_\emptyset$, $n_{k+1}=n_{(n_1, \dotsc ,n_k)}$ for $k\ge 1$. Then clearly the sequence $\{\vc_{(n_1)}, \dotsc ,\vc_{(n_1, \dotsc ,n_m)}, \dotsc \}$ witnesses the $\bs$-Hurewicz property of $X$. Since this is a sequence of moves of TWO in the $\bs$-Hurewicz game on $X$, $F$ is not a winning strategy for ONE in the $\bs$-Hurewicz game in $X$.

\noindent$(2)\Rightarrow (1)$.
Suppose that $X$ does not have the $\bk^s$-Hurewicz property. Then there exists a sequence $\{\uc_n:n\in \nb\}$ of open $\bk^s$-covers of $X$ such that  any sequence $\{\vc_n:n\in \nb\}$ of finite sets with $\vc_n\subseteq \uc_n$, fails to witness the $\bk^s$-Hurewicz property.
Let us define a strategy $F$ for ONE in the $\bs$-Hurewicz game on $X$ as follows. The first move of ONE is $F(X)=\uc_1$. TWO responds by choosing a finite set $\vc_1\subseteq \uc_1$. In the $n$-th inning ONE's move is $F(\vc_1,\vc_2, \dotsc,\vc_{n-1})=\uc_n\setminus (\vc_1\cup \vc_2\cup \cdots \cup \vc_{n-1})$. Since $F$ is not a wining strategy for ONE, the sequence $\{\vc_n:n\in \nb\}$ of moves of TWO witnesses the $\bk^s$-Hurewicz property, contradicting our assumption.
\ep

Next we define a $\gamma_{\bk^s}$-set.
\begin{Def}
\label{D2}
Let $\bk$ be a bornology on $X$ with closed base. $X$ is said to be a $\gamma_{\bk^s}$-set if $X$ satisfies the selection principle $\S1(\oc_{\bk^s},\Gamma_{\bk^s})$.
\end{Def}

By \cite[Theorem 2.8]{cmk}, $X$ is a $\gamma_{\bk^s}$-set if and only if every open $\bk^s$-cover of $X$ contains a countable set which is a $\gamma_{\bk^s}$-cover of $X$. Thus $X$ is a $\gamma_{\bk^s}$-set if and only if $X$ satisfies $\binom{\obs}{\gbs}$.

The following result can be easily verified.
\begin{Lemma}
\label{L4}
Every $\gamma_{\bk^s}$-set has the $\bk^s$-Hurewicz property.
\end{Lemma}


Next we introduce the following definition.
\begin{Def}
\label{D3}
Let $\bk$ be a bornology on $X$ with closed base. An open cover $\uc$ of $X$ is said to be $\bk^s$-groupable if it can be expressed as a union of countably many finite pairwise disjoint sets $\uc_n$ such that for each $B\in \bk$ there exist a $n_0\in \nb$ and a sequence $\{\delta_n:n\ge  n_0\}$ of positive real numbers with $B^{\delta_n}\subseteq U$ for some $U\in \uc_n$ for all $n\ge  n_0$.
\end{Def}
The collection of all $\bk^s$-groupable covers of $X$ is denoted by $\oc_{\bk^s}^{gp}$.\\

In the following result we show that $\bs$-Hurewicz property ensures the $\bs$-groupability of every (countable) open $\bs$-cover.
\begin{Lemma}
\label{L9}
Let $\bk$ be a bornology on $X$ with closed base. If $X$ has the $\bk^s$-Hurewicz property, then every (countable) open-$\bk^s$-cover of $X$ is $\bk^s$-groupable.
\end{Lemma}
\bp
Let $\uc$ be an open $\bk^s$-cover of $X$. By Proposition \ref{T1}, if we remove finitely many elements from $\uc$, it still remains an open $\bk^s$-cover.  Consider the strategy $\sigma$ in the $\bs$-Hurewicz game on $X$. Suppose that the first move of ONE is $\sigma(\emptyset)=\uc$. TWO responds with a finite set $\vc_1\subseteq \uc$. The second move of ONE is $\sigma(\vc_1)=\uc\setminus \vc_1$. TWO responds with a finite set $\vc_2\subseteq \sigma(\vc_1)$. Continuing, in the $n$-th inning, let the $n$-th move of ONE be $\sigma(\vc_1,\vc_2, \dotsc ,\vc_{n-1})=\uc\setminus (\vc_1\cup \vc_2\cup \cdots \cup \vc_{n-1})$.
By Theorem \ref{T8}, the play $\sigma(\emptyset), \vc_1, \sigma(\vc_1), \vc_2, \dotsc ,\sigma(\vc_1, \dotsc ,\vc_{n-1}), \vc_n, \dotsc$ is lost by ONE. So for each $B\in \bk$ there exist a  $n_0\in \nb$ and a sequence $\{\delta_n:n\ge n_0\}$ of positive real numbers such that $B^{\delta_n}\subseteq U$ for $U\in \vc_n$ for all $n\ge  n_0$ and hence $\cup_{n\in \nb}\vc_n$ is an open $\bk^s$-cover of $X$. Now $\{\vc_n:n\in \nb\}$ is so constructed that $\vc_n$'s are pairwise disjoint and finite. So $\{\vc_n:n\in \nb\}$ witnesses the groupability of $\cup_{n\in \nb}\vc_n$. Since $\uc$ is countable, the elements of $\uc\setminus \cup_{n\in \nb}\vc_n$ can be distributed among $\vc_n$'s so that $\{\vc_n:n\in \nb\}$ witnesses the $\bs$-groupability of  $\uc$. Hence $\uc$ is $\bk^s$-groupable cover of $X$.
\ep

The next result shows that $\Sf$-type selection hypothesis suffices to classify the $\bs$-Hurewicz property.
\begin{Th}
\label{T10}
Let $\bk$ be a bornology on $X$ with closed base. The following statements are equivalent.\\
\noindent$(1)$ $X$ has $\bk^s$-Hurewicz property.\\
\noindent$(2)$ $X$ satisfies $\Sf(\oc_{\bk^s},\oc_{\bk^s}^{gp})$.\\
\noindent$(3)$ ONE has no wining strategy in the game $\Gf(\oc_{\bk^s},\oc_{\bk^s}^{gp})$.
\end{Th}

\bp

\noindent$(2)\Rightarrow (1)$.
Let  $\{\uc_n:n\in \nb\}$ be a sequence of open $\bk^s$-covers of $X$. Enumerate $\uc_n$ bijectively as $\uc_n=\{U_k^n:k\in \nb\}$. Consider $\vc_n=\{U_{m_1}^1\cap \dotsc \cap U_{m_n}^n:n < m_1< \cdots < m_n\}$. By Lemma \ref{L2}, each $\vc_n$ is an open $\bk^s$-cover of $X$. Now apply $\Sf(\oc_{\bk^s},\oc_{\bk^s}^{gp})$ to the sequence $\{\vc_n:n\in \nb\}$, to obtain a sequence $\{\wc_n:n\in\nb\}$ of finite pairwise disjoint sets such that $\wc_n\subseteq \vc_n$ and $\cup_{n\in \nb} \wc_n$ is a $\bk^s$-groupable open cover of $X$. Let $\{\yc_n:n\in\nb\}$ witness the $\bs$-groupability of $\cup_{n\in \nb} \wc_n$. Clearly $\yc_n$'s are pairwise disjoint sets satisfying $\cup_{n\in \nb} \wc_n=\cup_{n\in \nb} \yc_n$.

Choose $n_1>1$ so large that $\yc_{n_1}\subseteq \cup_{j>1} \wc_j$ and let $\ac_1$ be the set of all $U_m^1$ that appear as $1$st component in the representation of elements of $\yc_{n_1}$. Choose $n_2 > n_1$  large enough so that $\yc_{n_2}\subseteq \cup_{j>2} \wc_j$ and let $\ac_2$ be set of all $U_m^2$ that appear as $2$nd component in the representation of elements of $\yc_{n_2}$. Continuing in this way at the $k$-th step we choose $\ac_k\subseteq \uc_k$. So we obtain a sequence $\{\ac_k:k\in \nb\}$ of finite sets such that $\ac_k\subseteq \uc_k$ for each $k$. For $B\in \bk$ there exist a $n_0\in \nb$ and a sequence $\{\delta_n:n\ge  n_0\}$ of positive real numbers such that $B^{\delta_n}\subseteq Y$ for all $n\ge  n_0$ for some $Y\in \yc_n$. Choose a $k_0$ such that $k\ge k_0$ implies $n_{k_0}>n_0$. Clearly $B^{\delta_{n_k}}\subseteq Y$ for some  $Y\in \yc_{n_k}$ for all $k\ge  k_0$. Define $\delta_k=\delta_{n_k}$ for each $k$. By the definition of $\ac_k$, we have $B^{\delta_k}\subseteq Y\subseteq U$ for some $U\in \ac_k$ for all $k\ge  k_0$. Hence $\{\ac_k:k\in \nb\}$ witnesses the $\bk^s$-Hurewicz property of $X$.\\

\noindent$(1)\Rightarrow (3)$.
Let $\tau$ be a strategy for ONE in $\Gf (\oc_{\bk^s},\oc_{\bk^s}^{gp})$.
Define a strategy $\sigma$ for ONE in the $\bs$-Hurewicz game as follows. Let $\uc$ be an open $\bk^s$-cover of $X$. Suppose that the first move of ONE is $\sigma(\emptyset)=\tau(\emptyset)=\uc$. TWO responds with a finite set $\vc_1\subset \uc$ (in the $\bs$-Hurewicz game). Then the second move of ONE is $\sigma(\vc_1)=\tau(\vc_1)\setminus \vc_1$. Let TWO respond with a finite set $\vc_2\subseteq \sigma(\vc_1)$. The $n$-th move of ONE is $\sigma(\vc_1, \dotsc, \vc_{n-1})=\tau(\vc_1, \dotsc, \vc_{n-1})\setminus (\vc_1\cup \dotsc \vc_{n-1})$ which is an open $\bs$-cover by Proposition \ref{P1}. TWO responds with a finite set  $\vc_n\subseteq \sigma(\vc_1, \dotsc, \vc_{n-1})$.  Continuing in this way we obtain the legitimate strategy $\sigma$ for ONE in the $\bs$-Hurewicz game on $X$. Since $X$ has $\bk^s$-Hurewicz property, the play $\sigma(\emptyset), \vc_1, \sigma(\vc_1,\vc_2),\vc_2, \dotsc$ is lost by ONE. Thus for each $B\in \bk$ there exist a $n_0\in \nb$ and a sequence $\{\delta_n:n\ge  n_0\}$  such that $B^{\delta_n}\subseteq V$ for all $n\ge  n_0$ for some $V\in \vc_n$. Clearly $\cup_{n\in \nb} \vc_n$ is a $\bk^s$-groupable open cover of $X$ as $\vc_n$'s are pairwise disjoint finite sets. Now for the strategy $\tau$, the play $\tau(\emptyset), \vc_1, \tau(\vc_1),\vc_2, \tau(\vc_1,\vc_2),\vc_3, \dotsc$ is legitimate in $\Gf(\oc_{\bk^s},\oc_{\bk^s}^{gp})$. As $\{\vc_n:n\in \nb\}$ is a sequence of moves by TWO in $\Gf(\oc_{\bk^s},\oc_{\bk^s}^{gp})$ and $\cup_{n\in \nb} \vc_n$ is a $\bk^s$-groupable open cover of $X$. Hence $\tau$ is not a wining strategy for ONE in $\Gf(\oc_{\bk^s},\oc_{\bk^s}^{gp})$.\\

\noindent$(3)\Rightarrow (2)$.
Let $X$ do not satisfy $\Sf(\oc_{\bk^s},\oc_{\bk^s}^{gp})$. Then there exists a sequence $\{\uc_n:n\in \nb\}$ of open $\bk^s$-covers of $X$ such that for any sequence $\{\vc_n:n\in \nb\}$, where each $\vc_n$ is a finite subset of $\uc_n$, we have $\cup_{n\in\nb} \vc_n\notin \oc_{\bk^s}^{gp}$.
Let us define a strategy $F$ for ONE in the game $\Gf(\oc_{\bk^s},\oc_{\bk^s}^{gp})$ in $X$. The first move of ONE is $F(X)=\uc_1$. TWO responds by choosing a finite set $\vc_1\subseteq \uc_1$. In the $n$-th inning ONE's move is $F(\vc_1,\vc_2, \dotsc,\vc_{n-1})=\uc_n$ and TWO responds by choosing  a finite set $\vc_n\subseteq \uc_n$. Since $F$ is not a wining strategy for ONE, we must have $\cup_{n\in\nb} \vc_n\in \oc_{\bk^s}^{gp}$, which contradicts our assumption. Hence $X$ satisfies $\Sf(\oc_{\bk^s},\oc_{\bk^s}^{gp})$.
\ep

\section{Results in function spaces}

This section is devoted to function spaces which naturally arise for metric spaces. Let $\bk$ be a bornology on $(X,d)$ with closed base and let $(Y,\rho)$ be another metric space. For $f\in C(X,Y)$, the neighbourhood of $f$ with respect to the topology $\tau_{\bk}^s$ of strong uniform convergence is denoted by
\[ [B,\varepsilon ]^s(f)=\{g\in C(X,Y):\exists \delta>0,\rho(f(x),g(x))<\varepsilon , \forall x\in B^\delta\}, \]
for $B\in\bk,\,\varepsilon >0$ (\!\cite{bl,cmh}).

The symbol $\underline{0}$ denotes the  zero function in $(C(X), \tau_{\bk}^s)$. The space $(C(X), \tau_{\bk}^s)$ is homogeneous and so it is enough to concentrate at the point $\underline{0}$ when dealing with local properties of this function space.
The following Lemma from \cite{cmk} will be used throughout this section.
\begin{Lemma}{(\cite[Lemma 2.2]{cmk})}
\label{koc}
Let  $\bk$ be a bornology on the metric space $(X,d)$. Consider the following statements.

\noindent $(a)$ Let $\uc$ be an open $\bs$-cover of $X$. If $A = \{f \in  C(X): \exists U \in\uc, f(x) = 1 \,\,\text{for all}\,\, x \in X \setminus U\}$. Then $\underline{0} \in \overline{A} \setminus A$ in $(C(X), \tau_\bk^s)$.\\
\noindent $(b)$ Let $A \subseteq (C(X), \tau_\bk^s)$ and let $\uc = \{f^{-1}(-\frac{1}{n},\frac{1}{n}): f \in A\}$, where $n\in\nb$. If $\underline{0} \in \overline{A}$ and $X\notin \uc$, then $\uc$ is an  open $\bs$-cover of $X$.
\end{Lemma}

\subsection{Certain applications of $\bs$-covers}

We start with a basic observation about $\bs$-covers.
\begin{Th}
\label{T13}
Let $\bk$ be a bornology on $X$ with closed base. Let $\uc$ be a collection of open subsets in $X$. The following statements are equivalent.\\
\noindent$(1)$ $\uc$ is an open $\bk^s$-cover of $X$.\\
\noindent$(2)$ For each $U\in \uc$ there is a closed set $C(U)\subseteq U$ such that $\{C(U):U\in \uc\}$ is a $\bk^s$-cover of $X$.
\end{Th}

\bp $(1)\Rightarrow (2)$
Let $\uc$ be an open $\bk^s$-cover of $X$. Consider the set $A=\{f\in C(X):\exists U\in \uc, f(X\setminus U)=\{1\}\}$.
By Lemma \ref{koc}, $A\in \Omega_{\underline{0}}$.
For $U\in \uc$, if $f\in A$ such that $f(X\setminus U)=\{1\}$, then $f^{-1}([-\frac{1}{2},\frac{1}{2}])\subseteq U$. Now take the closed set $C(U)$ to be $f^{-1}[-\frac{1}{2},\frac{1}{2}]$. Consider the collection $\vc=\{C(U):U\in \uc\}=\{f^{-1}[-\frac{1}{2},\frac{1}{2}]:f\in A\}$. We shall show that $\vc$ is a $\bk^s$-cover of $X$. For $B\in \bk$, consider the neighbourhood $[B,\frac{1}{2}]^s(\underline{0})$. Then $f\in [B,\frac{1}{2}]^s(\underline{0})\cap A\neq \emptyset$ i.e. there exists a $\delta>0$ such that $|f(x)|<\frac{1}{2}$ for all $x\in B^\delta$, i.e. $B^\delta\subseteq f^{-1}(-\frac{1}{2},\frac{1}{2})\subseteq f^{-1}[-\frac{1}{2},\frac{1}{2}]$. This shows that $\vc$ is a $\bk^s$-cover of $X$.

$(2)\Rightarrow (1)$ Can be similarly proved.
\ep


In the next few results we investigate how various properties of $C(X)$ can be characterized in terms of $\bs$-covers of $X$. We start with the property of countable $T$-tightness of $C(X)$ whose characterization is established in the next result.
\begin{Th}
\label{T12}
Let $\bk$ a bornology on $X$ with closed base. The following statements are equivalent.\\
\noindent$(1)$ $\left(C(X),\tau_\bk^s\right)$ has countable $T$-tightness.\\
\noindent$(2)$ For each uncountable regular cardinal $\rho$ and each increasing sequence $\{\uc_\alpha:\alpha<\rho\}$ of families of open subsets of $X$ such that $\bigcup_{\alpha<\rho}\uc_\alpha$ is an open ${\bk^s}$-cover of $X$, there is a $\beta<\rho$ with $\uc_\beta$ being an open $\bk^s$-cover of $X$.
\end{Th}

\bp

\noindent$(1)\Rightarrow (2)$.
Consider an increasing sequence of families of open subsets $\{\uc_\alpha:\alpha<\rho\}$ of $X$ such that $\bigcup_{\alpha<\rho}\uc_\alpha$ is an open $\bk^s$-cover of $X$. For $B\in \bk$, there exist a $\delta>0$ and $U\in \bigcup_{\alpha<\rho} \uc_\alpha$ such that $B^{2\delta}\subseteq U$. Define a continuous function $f_{B,U}$ on $X$ such that $f_{B,U}(B^\delta)=\{0\}$ and $f_{B,U}(X\setminus U)=\{1\}$. For each $\alpha<\rho$, consider $A_\alpha=\left\{f_{B,U}:U\in \uc_\alpha, B\in \bk\right\}$.
Clearly $\{\overline{A_\alpha}:\alpha< \rho\}$ is an increasing sequence of closed subsets of $C(X)$. By $(1)$, the set $A=\bigcup_{\alpha<\rho}\overline{A_\alpha}$ is closed in $\left(C(X),\tau_\bk^s\right)$.
Since  $\bigcup_{\alpha<\rho}\uc_\alpha$ is an open $\bk^s$-cover of $X$, $\underline{0}\in \overline{\cup_{\alpha<\rho} A_\alpha}\setminus \cup_{\alpha<\rho} A_\alpha$, by Lemma \ref{koc}$(a)$. Clearly $\underline{0}\in A$ and so one can find a $\beta<\rho$ with $\underline{0}\in \overline{A_\beta}$. By Lemma \ref{koc}$(b)$, $\{f_{B,U}^{-1}(-1,1):f_{B,U}\in A_\beta\}$ is an open $\bs$-cover of $X$. Since $f_{B,U}^{-1}(-1,1)\subseteq U\in \uc_\beta$, $\uc_\beta$ is an open $\bk^s$-cover of $X$.

\noindent$(2)\Rightarrow (1)$.
Consider an increasing sequence  $\{A_\alpha:\alpha<\rho\}$ of closed subsets of $C(X)$. We shall show that $A=\cup_{\alpha<\rho}A_\alpha$ is closed in $C(X)$. Without any loss of generality assume $\underline{0}\in \overline{A}$. Consider the collection $\uc_{\alpha,n}=\{f^{-1}(-\frac{1}{n},\frac{1}{n}):f\in A_\alpha\}$ and $\uc_n=\bigcup_{\alpha<\rho}\uc_{\alpha,n}$. By Lemma \ref{koc}$(b)$, $\uc_n$ is an open $\bk^s$-cover of $X$ for each $n\in \nb$. Now by our assumption, for each $n\in \nb$ there is an open $\bk^s$-cover $\uc_{{\beta_n},n}\subset \uc_n$. Let $\beta_0=\sup\{\beta_n:n\in \nb\}$.
It can be easily seen that for each $n$, $\uc_{\beta_0,n}$ is an open $\bk^s$-cover of $X$.
Now we show that $\underline{0}\in A_{\beta_0}$. Let $B\in \bk$ and $\varepsilon>0$ be given. Choose  $n_0\in \nb$  such that $\frac{1}{n}<\varepsilon$ for all $n\ge  n_0$. Since $\uc_{\beta_0,n}=\{f^{-1}(-\frac{1}{n},\frac{1}{n}):f\in A_{\beta_0}\}$ is an open $\bk^s$-cover of $X$, there exist a $f\in A_{\beta_0}$ and a $\delta>0$ such that $B^\delta\subseteq f^{-1}(-\frac{1}{n},\frac{1}{n})\subseteq f^{-1}(-\varepsilon,\varepsilon)$ for $n\ge  n_0$. Thus $f\in [B,\varepsilon]^s(\underline{0})\cap A_{\beta_0}\neq \emptyset$ and so $\underline{0}\in \overline{A_{\beta_0}}=A_{\beta_0}$, i.e. $\underline{0}\in A$. This proves that $A$ is closed in $X$.
\ep

Using  techniques of the proof of  \cite[Theorem 12]{s} and using Lemma \ref{koc} one can also characterize weakly Fr\'{e}chet-Urysohn property as well.

\begin{Th}
\label{15}
Let $\bk$ be a bornology on $X$ with closed base. Then the following statements are equivalent.\\
\noindent$(1)$ $(C(X),\tau_{\bk}^s)$ is weakly Fr\'{e}chet-Urysohn.\\
\noindent$(2)$ Every open $\bk^s$-cover of $X$ is  $\bk^s$-groupable.
\end{Th}

In the next result we obtain a characterization of a countable $\pi$-network at $\underline 0$ of finite subsets of  $A\subset \left(C(X),\tau_\bk^s\right)$.
\begin{Th}
\label{T14m}
Let $\bk$ be a bornology on $X$ with closed base. The following statements are equivalent.\\
\noindent$(1)$ If $A$ is subset of $\left(C(X),\tau_\bk^s\right)$ and if $\underline{0}\in \overline{A}\setminus A$, then there is a countable $\pi$-network at $\underline{0}$ of finite subsets of $A$.\\
\noindent$(2)$ If $\uc$ is an open $\bk^s$-cover of $X$, then there is a sequence $\{\uc_n:n\in \nb\}$ of finite subfamilies of $\uc$ such that $\{\bigcap\uc_n:n\in \nb\}$ is an open $\bk^s$-cover of $X$.
\end{Th}

\bp
\noindent$(1)\Rightarrow (2)$.
Let $\uc$ be an open $\bk^s$-cover of $X$. Consider the set $A=\{f_U\in C(X):\exists U\in \uc, f_U(X\setminus U)=\{1\}\}$. By Lemma \ref{koc}$(a)$, $A\in \Omega_{\underline{0}}$. By $(1)$ there is a sequence $\{A_n:n\in \nb\}$ of finite subsets of $A$ such that $\{A_n:n\in \nb\}$ is a $\pi$-network at $\underline{0}$. Now choose a finite subset $\uc_n$ of $\uc$ such that $A_n=\{f_U:U\in \uc_n\}$. We show that $\{\bigcap\uc_n:n\in \nb\}$ is an open $\bk^s$-cover of $X$.
Let $B\in \bk$ and consider the neighbourhood $[B,1]^s(\underline{0})$ of $\underline{0}$. Then $A_n\subset [B,1]^s(\underline{0})$ for some $n\in \nb$ as $\{A_n:n\in \nb\}$ is a $\pi$-network at $\underline{0}$. Clearly $f_U\in [B,1]^s(\underline{0})$ for all $U\in \uc_n$. If $\uc_n=\{U_n^1, U_n^2, \cdots, U_n^{r_n}\}$, then for each $k=1,2, \dotsc,r_n$, there exists a $\delta_k>0$ such that $|f_{U_n^k}(x)|<1$ for all $x\in B^{\delta_k}$, i.e. $B^{\delta_k}\subseteq U_n^k$. Consequently there exists a  $\delta>0$ such that $B^\delta\subseteq \bigcap \uc_n$. Hence $\{\bigcap \uc_n:n\in \nb\}$ is an open $\bk^s$-cover of $X$.

$(2)\Rightarrow (1)$.
Let $A\subset C(X)$ be such that $\underline{0}\in \overline{A}\setminus A$.
For each $n\in \nb$, consider the collection $\uc_n=\{f^{-1}(-\frac{1}{n},\frac{1}{n}):f\in A\}$. Then by Lemma \ref{koc}$(b)$, $\uc_n$ is an open $\bk^s$-cover of $X$. \\
If there exists a sequence $n_1<n_2< \cdots $ of positive integers such that $f_k^{-1}(-\frac{1}{n_k},\frac{1}{n_k})=X$, then the sequence $\{f_k\}_{k\in \nb}$ converges to $\underline{0}$. Choose $A_1=\{f_i:i<n_1\}$ and for each $k>1$ choose $A_k=\{f_i:n_{k-1}\leq i<n_k\}$. Consider the family $\ac=\{A_k:k\in \nb\}$. For $B\in \bk$ and $\epsilon>0$, consider the neighbourhood $[B,\epsilon]^s(\underline{0})$ of $\underline{0}$. Choose $k_0\in \nb$ such that $\frac{1}{n_k}< \epsilon$ for all $k\ge k_0$. Now for any $\delta>0$, $B^\delta\subseteq f_k^{-1}(-\frac{1}{n_k},\frac{1}{n_k})=X$ for all $k\in \nb$ and hence $f_k\in [B,\epsilon]^s(\underline{0})$ for all $k\ge k_0$, i.e. $A_{k+1}\subset [B,\epsilon]^s(\underline{0})$ for all $k\ge k_0$. Hence $\ac=\{A_k:k\in \nb\}$ is a $\pi$-network at $\underline{0}$.

Otherwise, we assume that there is a $n_0\in \nb$ such that $X\notin \uc_n$ for all $n\ge  n_0$. Now by $(2)$, for each $n\ge n_0$, there is a sequence $\{\vc_{n,m}:m\in \nb\}$ of finite subfamilies of $\uc_n$ such that $\{\bigcap\vc_{n,m}:m\in \nb\}$ is an open $\bk^s$-cover of $X$.
Consequently there is a sequence $\{A_{n,m}:n\ge n_0, m\in \nb\}$ of finite subsets of $A$ such that $\vc_{n,m}=\{f^{-1}(-\frac{1}{n},\frac{1}{n}):f\in A_{n,m}\}$ for each $n\ge  n_0$ and $m\in\nb$.
Consider $\ac=\{A_{n,m}:n\ge  n_0, m\in \nb\}$. Let $[B,\varepsilon]^s(\underline{0})$ be a neighbourhood of $(\underline{0})$, where $\epsilon>0$. Choose  $n_1\in \nb$ such that $\frac{1}{n}<\epsilon$ for all $n\ge n_1$. Fix  $n\ge \max\{n_0,n_1\}$. Then there exists a $\delta>0$ such that $B^\delta\subseteq \cap \vc_{n,m}$ for some $m\in \nb$, i.e. $B^\delta\subseteq f^{-1}(-\frac{1}{n},\frac{1}{n})$ for all $f\in A_{n,m}$ for some $m$, i.e. $f\in [B,\frac{1}{n}]^s(\underline{0})$ for all $f\in A_{n,m}$ for some $m$. So $A_{n,m}\subseteq [B,\frac{1}{n}]^s(\underline{0})\subseteq [B,\varepsilon]^s(\underline{0})$ for some $m$ and each $n\ge  \max\{n_0,n_1\}$. Hence $\ac= \{A_{n,m}: n\ge  n_0, m\in \nb\}$ is a $\pi$-network at $\underline{0}$.
\ep

Similarly we can prove a sufficient condition for $C(X)$ to be Pytkeev.
\begin{Th}
\label{T14}
Let $\bk$ be a bornology on $X$ with closed base. Consider the following statements.\\
\noindent$(1)$ If $\uc$ is an open $\bk^s$-cover of $X$, then there is a sequence $\{\uc_n:n\in \nb\}$ of countably infinite subfamilies of $\uc$ such that $\{\bigcap\uc_n:n\in \nb\}$ is an open $\bk^s$-cover of $X$.\\
\noindent$(2)$ $\left(C(X),\tau_\bk^s\right)$ is  Pytkeev.\\

Then $(1)\Rightarrow (2)$ holds.
\end{Th}

\begin{Rem}
We do not know whether the preceding condition is necessary for $C(X)$ to be Pytkeev. We leave it as an open problem.
\end{Rem}

\begin{Problem}
Determine whether the condition $(1)$ in Theorem \ref{T14} is necessary for $\left(C(X),\tau_\bk^s\right)$ to be Pytkeev.
\end{Problem}

\subsection{Game theoretic results in function spaces}


In the next two results we observe the connection between game theoretic relation on $X$ and selection hypothesis on $C(X)$.
First we state a result from \cite{cmk}. It was shown \cite[Theorem 2.3]{cmk} that
\emph{$(C(X),\tau_{\bk}^s)$ has countable strong fan tightness if and only if $X$ satisfies $\S1(\oc_{\bk^s},\oc_{\bk^s})$.
}
\begin{Th}
\label{T16}
Let $\bk$ be a bornology on $X$ with closed base. If ONE has no wining strategy in the game $\G1 (\oc_{\bk^s},\oc_{\bk^s}^{gp})$ on $X$, then $(C(X),\tau_\bk^s)$ satisfies $S_1(\Omega_{\underline{0}},\Omega_{\underline{0}}^{gp})$.
\end{Th}

\bp
Clearly if $X$ satisfies $\S1(\oc_{\bk^s},\oc_{\bk^s}^{gp})$ then $X$ also satisfies $\S1(\oc_{\bk^s},\oc_{\bk^s})$. By \cite[Theorem 2.3]{cmk}, $(C(X),\tau_{\bk}^s)$ has countable strong fan tightness, i.e. $(C(X),\tau_{\bk}^s)$ satisfies $\S1(\Omega_{\underline{0}},\Omega_{\underline{0}})$.
Now it suffices to show that each countable subset of $\Omega_{\underline{0}}$ is groupable. Let $A\in \Omega_{\underline{0}}$ be a countable subset of $C(X)$. By Lemma \ref{koc}, $\uc_1=\{f^{-1}(-1,1):f\in A\}$ is an open $\bk^s$-cover of $X$.
We define a strategy $\sigma$ for ONE in $\G1 (\oc_{\bk^s},\oc_{\bk^s}^{gp})$ as follows. Let the first move of ONE be $\sigma(\emptyset)=\uc_1$. TWO chooses $U_1\in \uc_1$. Take the corresponding function $f_1\in A$ such that $U_1=f_1^{-1}(-1,1)$. Now consider $A_1=A\setminus \{f_1\}$ with $\underline{0}\in \overline{A_1}$. Again by Lemma \ref{koc}, $\uc_2=\{f^{-1}(-\frac{1}{2},\frac{1}{2}):f\in A_1\}\setminus \{U_1\}$ is an open $\bk^s$-cover of $X$. Let the 2nd move of ONE be $\sigma(U_1)=\uc_2$. Take the corresponding function $f_2\in A_1$ such that $U_2=f_2^{-1}(-\frac{1}{2},\frac{1}{2})$. Continuing in this way we have the following.

\begin{enumerate}
\item[$(a)$] $\{\uc_n:n\in \nb\}$ is a sequence of open $\bk^s$-cover of $X$, where $\uc_n=\sigma(U_1,U_2, \dotsc ,U_{n-1})$.
\item[$(b)$] For each $n\in \nb$, $U_n$ is not in $\{U_1,U_2, \dotsc, U_{n-1}\}$.
\item[$(c)$] For each $n\in \nb$, $f_n\in A\setminus \{f_1,f_2, \dotsc, f_{n-1}\}$.
\item[$(d)$] For each $n\in \nb$, $U_n=f_n^{-1}(-\frac{1}{n},\frac{1}{n})$.
\end{enumerate}
Since $\sigma$ is not a wining strategy for ONE, the play $\uc_1,U_1 \dotsc, \uc_n,U_n, \dotsc $ is lost by ONE and consequently  $\vc=\{U_n:n\in \nb\}$ is a $\bk^s$-groupable cover of $X$. Therefore there is an increasing infinite sequence $n_1<n_2<  \cdots$ such that the sets $\hc_k=\{U_i:n_k\le  i\le  n_{k+1}\}$ (for $k=1,2, \dotsc$) are finite, pairwise disjoint and for every $B\in \bk$ there exist a $k_0\in \nb$ and  a sequence $\{\delta_k:k\ge  k_0\}$ of positive real numbers such that $B^{\delta_k}\subseteq U$ for $U\in \hc_k$ for all $k\ge  k_0$. Define $\mc_k=\{f_i:n_k\le  i<n_{k+1}\}$. Then $\mc_k$'s are finite, pairwise disjoint subsets of $A$. Clearly $A=\bigcup_{k\in \nb}\mc_k$. We claim that the sequence $\{\mc_k:k\in \nb\}$ witnesses the groupability of $A$.\\

To see this for $B\in \bk$, consider the neighbourhood $[B,\varepsilon ]^s(\underline{0})$ of $\underline{0}$. Choose a $k_1\in \nb$ such that $\frac{1}{k}<\varepsilon $ for all $k\ge  k_1$. Again for that $B$ there exist a $k_0\in \nb$ and a sequence $\{\delta_k:k\ge  k_0\}$ of positive real numbers such that $B^{\delta_k}\subseteq U$ for $U\in \hc_k$ for all $k\ge  k_0$, i.e. $B^{\delta_k}\subseteq f^{-1}(-\frac{1}{k}, \frac{1}{k})$. Observe that $|f(x)|<\frac{1}{k}<\varepsilon $ for all $x\in B^{\delta_k}$ and $k\ge  k_2$ where $k_2 \ge \max\{k_0,k_1\}$ i.e. $f\in \mc_k\cap [B,\varepsilon ]^s(\underline{0})\neq \emptyset$ for all $k\ge  k_2$. Hence every neighbourhood of $\underline{0}$ intersect all but finitely many $\mc_k$. So $A\in \Omega_{\underline{0}}^{gp}$. Hence $X$ satisfies $S_1(\Omega_{\underline{0}},\Omega_{\underline{0}}^{gp})$.
\ep

Using similar argument we can prove the following result.

\begin{Th}
\label{T17}
Let $\bk$ be a bornology on $X$ with closed base. If ONE has no wining strategy in the game $\Gf(\oc_{\bk^s},\oc_{\bk^s}^{gp})$, then $(C(X),\tau_\bk^s)$ has property $\Sf(\Omega_{\underline{0}},\Omega_{\underline{0}}^{gp})$.
\end{Th}

In \cite[Theorem 21]{coocvii} it was proved that
\emph{
$C_p(X)$ has countable fan tightness and Reznichenko's property if and only if $C_p(X)$ has the property $\Sf(\Omega_{\underline{0}}\Omega_{\underline{0}}^{gp})$.} Using the techniques of this proof and keeping in mind Theorem \ref{T10} and Theorem \ref{T16}, we can prove the next result.

\begin{Th}
\label{T18}
Let $\bk$ be  bornology on $X$ with closed base. The following statements are equivalent.\\
\noindent$(1)$ $X$ has the $\bk^s$-Hurewicz property.\\
\noindent$(2)$ $(C(X),\tau_{\bk}^s)$ has countable fan tightness and Reznichenko property.
\end{Th}

Let $\bk$ be a bornology on $X$ with closed base. The game $\G1 (\Omega_{\underline{0}},\Sigma_{\underline{0}})$ on $C(X)$ is defined as follows. Let $\psi$ be a strategy for ONE in the game $\G1 (\Omega_{\underline{0}},\Sigma_{\underline{0}})$. The first move of ONE is $\psi(\emptyset)=A_1 \in \Omega_{\underline{0}}$ and TWO responds by selecting $f_1\in A_1$. In the $n$-th inning ONE's move is $\psi(f_1,f_2, \dotsc ,f_{n-1})=A_{n}$ and TWO responds by choosing $f_n\in A_{n}$. TWO wins the play
$\psi(\emptyset), f_1,\psi(f_1), f_2, \dotsc ,\psi(f_1, \dotsc ,f_{n-1}), \dotsc$, if $\{f_n:n\in \nb\}\in \Sigma_{\underline{0}}$. Otherwise ONE wins.


Our final result concerns the above mentioned game and the idea of the proof follows that of \cite[Theorem 26]{roc}.

\begin{Th}
\label{T19}
Let $\bk$ be a bornology on $X$ with closed base. The following statements are equivalent.\\
\noindent$(1)$ TWO has a wining strategy in the game $\G1 (\Omega_{\underline{0}},\Sigma_{\underline{0}})$ on $(C(X),\tau_{\bk}^s)$.\\
\noindent$(2)$ TWO has a wining strategy in the game $\G1 (\oc_{\bk^s},\Gamma_{\bk^s})$ on $X$.
\end{Th}

\bp
\noindent$(1)\Rightarrow (2)$.
Let $\psi$ be a wining strategy for TWO in $\G1 (\Omega_{\underline{0}},\Sigma_{\underline{0}})$. We define a strategy $\sigma$ for ONE in the game $\G1 (\oc_{\bk^s},\Gamma_{\bk^s})$. ONE's first move in $\G1 (\oc_{\bk^s},\Gamma_{\bk^s})$ is $\sigma(\emptyset)=\uc_1\in \oc_{\bk^s}$. For each $B\in \bk$, there exist a $\delta>0$ and $U\in \uc_1$ such that $B^{2\delta}\subseteq U$. Take a continuous function $f_{B,U}$ on $X$ such that $f_{B,U}(B^{\delta})=\{0\}$ and $f_{B,U}(X\setminus U)=\{1\}$. Then the collection $A_1=\{f_{B,U}:B\in \bk, U\in \uc_1\}$ is in $\Omega_{\underline{0}}$ by Lemma \ref{koc}. ONE's first move in $\G1 (\Omega_{\underline{0}},\Sigma_{\underline{0}})$  is $\psi(\emptyset)=A_1$. TWO responds with $f_{B_1,U_1}\in A_1$. Now TWO's response in $\G1 (\oc_{\bk^s},\Gamma_{\bk^s})$ is $U_1$. In the $n$-th inning, the move of ONE in $\G1 (\oc_{\bk^s},\Gamma_{\bk^s})$ is $\sigma(U_1, \dotsc ,U_{n-1})=\uc_n$, then the corresponding move of ONE in
$\G1 (\Omega_{\underline{0}},\Sigma_{\underline{0}})$ is $A_n=\{f_{B,U}:B\in \bk, U\in \uc_n\}$. TWO responds with $f_{B_n,U_n}\in A_n$ in $\G1 (\Omega_{\underline{0}},\Sigma_{\underline{0}})$. Correspondingly in $\G1 (\oc_{\bk^s},\Gamma_{\bk^s})$, TWO's response is $U_n$.

The play in $\G1 (\oc_{\bk^s},\Gamma_{\bk^s})$ is
\[\sigma(\emptyset), U_1, \sigma(U_1), U_2, \dotsc ,\sigma(U_1, \dotsc ,U_{n-1}), U_n \dotsc\] and
the play in $\G1 (\Omega_{\underline{0}},\Sigma_{\underline{0}})$ is
\[\psi(\emptyset), f_1,\psi(f_1), f_2, \dotsc ,\psi(f_1, \dotsc ,f_{n-1}), f_n\dotsc,\]     where $f_n= f_{B_n,U_n}$.
Now by $(1)$, $\psi$ is a wining strategy for TWO in the play in $\G1 (\Omega_{\underline{0}},\Sigma_{\underline{0}})$. So $\{f_n:n\in \nb\}\in \Sigma_{\underline{0}}$. We want to show that $\{U_n:n\in \nb\}$ is a $\gamma_{\bk^s}$-cover of $X$. Let $B\in \bk$ and $\varepsilon =1$. For the neighbourhood $[B,1]^s(\underline{0})$ there exists a $n_0\in \nb$ such that $f_n\in [B,1]^s(\underline{0})$ for all $n\ge  n_0$, i.e. for each $n$, there exists a $\delta_n>0$ such that $|f_n(x)|<1$ for all $x\in B^{\delta_n}$, i.e. $B^{\delta_n}\subseteq U_n$ for all $n\ge  n_0$. So $\{U_n:n\in \nb\}\in \Gamma_{\bk^s}$. Hence TWO has a wining strategy in the game $\G1 (\oc_{\bk^s},\Gamma_{\bk^s})$.\\

\noindent$(2)\Rightarrow (1)$.
Let $\sigma$ be a wining strategy for TWO in the game $\G1 (\oc_{\bk^s},\Gamma_{\bk^s})$. For each $n$ let $I_n=(-\frac{1}{n+1},\frac{1}{n+1})$.
Now define a strategy $\psi$ for ONE in the $\G1 (\Omega_{\underline{0}},\Sigma_{\underline{0}})$. In the $n$-th inning ONE's move in $\G1 (\Omega_{\underline{0}},\Sigma_{\underline{0}})$ is $\psi(f_1, \dotsc ,f_{n-1})=A_n\in \Omega_{\underline{0}}$. Then $n$-th move of ONE in $\G1 (\oc_{\bk^s},\Gamma_{\bk^s})$ is $\uc(A_n)=\{f^{-1}(I_n):f\in A_n\}$, since $\uc(A_n)$ is an open $\bs$-cover by Lemma \ref{koc}. TWO responds with $U_n\in \uc(A_n)$. TWO's move in $\G1 (\Omega_{\underline{0}},\Sigma_{\underline{0}})$ is $f_n$, where $U_n=f_n^{-1}(I_n)$.
Now the play in $\G1 (\Omega_{\underline{0}},\Sigma_{\underline{0}})$ is
\[A_1, f_1, A_2, f_2, \dotsc ,A_n, f_n, \dotsc\]
and
the play in $\G1 (\oc_{\bk^s},\Gamma_{\bk^s})$ is
\[\uc(A_1),U_1, \uc(A_2), U_2, \dotsc ,\uc(A_n), U_n, \dotsc.\]
If for $A_n\in \Omega_{\underline{0}}$, $X\in \uc(A_n)$ for infinitely many $n$, then the conclusion is trivial. So we assume that $X\notin \uc(A_n)$ for all $n\ge  n_0$, $n_0\in \nb$. Since $\sigma$ is a wining strategy for TWO in $\G1 (\oc_{\bk^s},\Gamma_{\bk^s})$,  we have $\{U_n:n\in \nb\}\in \Gamma_{\bk^s}$. We show that $\{f_n:n\in \nb\}\in \Sigma_{\underline{0}}$.
Let $B\in \bk$ and choose a neighbourhood $[B,\varepsilon ]^s(\underline{0})$ of $\underline{0}$. Choose $n_1\in \nb$ such that $\frac{1}{n+1}<\varepsilon $ for $n\ge  n_1$. Now there is a $n_2\in \nb$ and a sequence $\{\delta_n:n\geqslant n_2\}$ such that $B^{\delta_n}\subseteq U_n$ for all $n\ge  n_2$. Since $U_n=f_n^{-1}(I_n)=f_n^{-1}(-\frac{1}{n+1},\frac{1}{n+1})$, we have $|f_n(x)|<\frac{1}{n+1}<\varepsilon $ for all $n\ge  n_1$ and $x\in U_n$. If $N=\max\{n_1,n_2\}$, then $|f_n(x)|<\varepsilon $ for all $x\in B^{\delta_n}$ and $n\ge  N$. Consequently $f_n\in [B,\varepsilon ]^s(\underline{0})$ for all $n\ge  N$, i.e. $\{f_n:n\in \nb\}\in \Sigma_{\underline{0}}$. Hence TWO has a wining strategy in the game $\G1 (\Omega_{\underline{0}},\Sigma_{\underline{0}})$.
\ep

{}
\end{document}